\documentclass{article}%
\usepackage{amsfonts}
\usepackage{graphicx}
\usepackage{amsmath}%
\usepackage{amssymb}
\providecommand{\U}[1]{\protect\rule{.1in}{.1in}}
\providecommand{\U}[1]{\protect\rule{.1in}{.1in}}
\providecommand{\U}[1]{\protect\rule{.1in}{.1in}}
\providecommand{\U}[1]{\protect\rule{.1in}{.1in}}
\providecommand{\U}[1]{\protect\rule{.1in}{.1in}}
\providecommand{\U}[1]{\protect\rule{.1in}{.1in}}
\newtheorem{theorem}{Theorem}

\newtheorem{conjecture}[theorem]{Conjecture}
\newtheorem{corollary}[theorem]{Corollary}

\newtheorem{lemma}[theorem]{Lemma}

\newtheorem{proposition}[theorem]{Proposition}

\newenvironment{proof}[1][Proof]{\textbf{#1.} }{\ \rule{0.5em}{0.5em}}
\begin{document}

\title{Microlinearity in Fr\"{o}licher Spaces\\-Beyond the Regnant Philosophy of Manifolds-}
\author{Hirokazu Nishimura\\Institute of Mathematics, University of Tsukuba\\Tsukuba, Ibaraki, 305-8571, Japan}
\maketitle

\begin{abstract}
Fr\"{o}licher spaces and smooth mappings form a cartesian closed category. It
was shown in our previous paper [Far East Journal of Mathematical Sciences, 35
(2009), 211-223] that its full subcategory of Weil exponentiable Fr\"{o}licher
spaces is cartesian closed. By emancipating microlinearity from within a
well-adapted model of synthetic differential geometry to Fr\"{o}licher spaces,
we get the notion of microlinearity for Fr\"{o}licher spaces. It is shown in
this paper that its full subcategory of Weil exponentiable and microlinear
Fr\"{o}licher spaces is cartesian closed. The canonical embedding of Weil
exponentiable Fr\"{o}licher spaces into the Cahiers topos is shown to preserve
microlinearity besides finite products and exponentiation.

\end{abstract}

\section{Introduction}

Differential geometry of finite-dimensional smooth manifolds has been
generalized by many authors to the infinite-dimensional case by replacing
finite-dimensional vector spaces by Hilbert spaces, Banach spaces, Fr\'{e}chet
spaces or, more generally, convenient vector spaces as the local prototype. We
know well that the category of smooth manifolds of any kind, whether
finite-dimensional or infinite-dimensional, is not cartesian closed, while
Fr\"{o}licher spaces, introduced by Fr\"{o}licher and others (cf. \cite{fro1},
\cite{fro2} and \cite{fro}), do form a cartesian closed category. We are
strongly biased in favor of our central dogma that \textit{the basic objects
under study of infinite-dimensional differential geometry should form a
cartesian closed category}. It seems that Fr\"{o}licher and his followers do
not know what a kind of Fr\"{o}licher space, besides convenient vector spaces,
should become the basic object of research for infinite-dimensional
differential geometry. The category of Fr\"{o}licher spaces and smooth
mappings should be restricted \textit{adequately} to a cartesian closed subcategory.

Synthetic differential geometry is differential geometry with a cornucopia of
nilpotent infinitesimals. For a standard textbook on synthetic differential
geometry the reader is referred to \cite{kock}, whose Chapter III is devoted
to its model theory. Roughly speaking, a space of nilpotent infinitesimals of
some kind, which exists only within an imaginary world, corresponds to a Weil
algebra, which is an entity of the real world. The central object of study in
synthetic differential geometry is \textit{microlinear} spaces. Although the
notion of a manifold (=a pasting of copies of a certain linear space) is
defined on the local level, the notion of microlinearity is defined absolutely
on the genuinely infinitesimal level. What we should do so as to get an
adequately restricted cartesian closed category of Fr\"{o}licher spaces is to
emancipate microlinearity from within a well-adapted model of synthetic
differential geometry. In other words, we should externalize the notion of
microlinearity for Fr\"{o}licher spaces, which is the princcipal objective in
this paper.

Although nilpotent infinitesimals exist only within a well-adapted model of
synthetic differential geometry, the notion of Weil functor was formulated for
finite-dimensional manifolds (cf. Section 35 of \cite{kolar}) and for
infinite-dimensional manifolds (cf. Section 31 of \cite{kri}). It was
generalized to that for Fr\"{o}licher spaces in our previous paper
\cite{nishi}, which paved the way to microlinearity for Fr\"{o}licher spaces.
This is the first step towards microlinearity for Fr\"{o}licher spaces.
Therein all Fr\"{o}licher spaces which believe in fantasy that all Weil
functors are really exponentiations by some adequate infinitesimal objects in
imagination form a cartesian closed category. This is the second step towards
microlinearity for Fr\"{o}licher spaces. We will introduce the notion of
''transversal limit diagram of Fr\"{o}licher spaces''\ after the manner of
that of ''transversal pullback''\ in Section \ref{2.5}, which is a familiar
token in the arena of synthetic differential geometry. This is the third and
final step towards microlinearity for Fr\"{o}licher spaces. Just as
microlinearity is closed under arbitrary limits within a well-adapted model of
synthetic differential geometry, microlinearity for Fr\"{o}licher spaces is
closed under arbitrary transversal limits. We will introduce the central
notion of microlinearity in Section \ref{3}, where it is to be shown that Weil
exponentiable and microlinear Fr\"{o}licher spaces, together with smooth
mappings among them, form a cartesian closed category. In Section \ref{4}\ we
will demonstrate that our canonical embedding of the category of Fr\"{o}licher
spaces and smooth mappings into the Cahiers topos preserves microlinearity.
Therein our hasty discussions in Section 5 of \cite{nishi} will also be elaborated.

\section{Preliminaries\label{2}}

\subsection{Weil Prolongation}

In our previous paper \cite{nishi} we have discussed how to assign, to each
pair $(X,W)$\ of a Fr\"{o}licher space $X$ and a Weil algebra $W$,\ another
Fr\"{o}licher space $X\otimes W$,\ called the \textit{Weil prolongation of}
$X$ \textit{with respect to} $W$, which naturally extends to a bifunctor
$\mathbf{FS}\times\mathbf{W\rightarrow FS}$, where $\mathbf{FS}$\ is the
category of Fr\"{o}licher spaces and smooth mappings, and $\mathbf{W}$\ is the
category of Weil algebras. We have shown in \cite{nishi}\ that

\begin{theorem}
\label{t2.1}The functor $\cdot\otimes W:\mathbf{FS\rightarrow FS}$ is
product-preserving for any Weil algebra $W$.
\end{theorem}

\subsection{Weil Exponentiability}

A Fr\"{o}licher space $X$ is called \textit{Weil exponentiable }if
\begin{equation}
(X\otimes(W_{1}\otimes_{\infty}W_{2}))^{Y}=(X\otimes W_{1})^{Y}\otimes
W_{2}\label{2.1}%
\end{equation}
holds naturally for any Fr\"{o}licher space $Y$ and any Weil algebras $W_{1}$
and $W_{2}$. If $Y=1$, then (\ref{2.1}) degenerates into
\begin{equation}
X\otimes(W_{1}\otimes_{\infty}W_{2})=(X\otimes W_{1})\otimes W_{2}\label{2.2}%
\end{equation}
If $W_{1}=\mathbb{R}$, then (\ref{2.1}) degenerates into
\begin{equation}
(X\otimes W_{2})^{Y}=X^{Y}\otimes W_{2}\label{2.3}%
\end{equation}

The following propositions and theorem have been established in our previous
paper \cite{nishi}.

\begin{proposition}
\label{t2.2.0}Convenient vector spaces are Weil exponentiable.
\end{proposition}

\begin{corollary}
$C^{\infty}$-manifolds in the sense of \cite{kri} (cf. Section 27) are Weil exponentiable.
\end{corollary}

\begin{proposition}
\label{t2.2.1}If $X$ is a Weil exponentiable Fr\"{o}licher space, then so is
$X\otimes W$ for any Weil algebra $W$.
\end{proposition}

\begin{proposition}
\label{t2.2.2}If $X$ and $Y$ are Weil exponentiable Fr\"{o}licher spaces, then
so is $X\times Y$.
\end{proposition}

\begin{proposition}
\label{t2.2.3}If $X$ is a Weil exponentiable Fr\"{o}licher space, then so is
$X^{Y}$ for any Fr\"{o}licher space $Y$.
\end{proposition}

\begin{theorem}
\label{t2.2.4}Weil exponentiable Fr\"{o}licher spaces, together with smooth
mappings among them, form a Cartesian closed subcategory $\mathbf{FS}%
_{\mathbf{WE}}$\ of the category $\mathbf{FS}$.
\end{theorem}

\section{\label{2.5}Transversal Limit Diagrams}

Generally speaking, limits in the category $\mathbf{FS}$ are bamboozling. The
notion of limit in $\mathbf{FS}$ should be elaborated geometrically.

A finite cone $\mathcal{D}$ in $\mathbf{FS}$ is called a \textit{transversal
limit} \textit{diagram} providing that $\mathcal{D}\otimes W$ is a limit
diagram in $\mathbf{FS}$ for any Weil algebra $W$, where the diagram
$\mathcal{D}\otimes W$ is obtained from $\mathcal{D}$ by putting $\otimes W$
to the right of every object and every morphism in $\mathcal{D}$. By taking
$W=\mathbb{R}$, we see that a transversal limit diagram is always a limit
diagram. The limit of a finite diagram of Fr\"{o}licher spaces is said to be
\textit{transversal} providing that its limit diagram is a transversal limit diagram.

\begin{lemma}
\label{t2.3.1}If $\mathcal{D}$ is a transversal limit diagram whose objects
are all Weil exponentiable, then $\mathcal{D}^{X}$ is also a transversal limit
diagram for any Fr\"{o}licher space $X$, where $\mathcal{D}^{X}$ is obtained
from $\mathcal{D}$ by putting $X$ as the exponential over every object and
every morphism in $\mathcal{D}$.
\end{lemma}

\begin{proof}
Since the functor $\cdot^{X}:\mathbf{FS}\rightarrow\mathbf{FS}$ preserves
limits, we have
\[
\mathcal{D}^{X}\otimes W=(\mathcal{D}\otimes W)^{X}%
\]
for any Weil algebra $W$, so that we have the desired result.
\end{proof}

\begin{lemma}
\label{t2.3.2}If $\mathcal{D}$ is a transversal limit diagram whose objects
are all Weil exponentiable, then $\mathcal{D}\otimes W$ is also a transversal
limit diagram for any Weil algebra $W$.
\end{lemma}

\begin{proof}
Since the functor $W\otimes_{\infty}\cdot:\mathbf{W\rightarrow W}$ preserves
finite limits, we have
\[
(\mathcal{D}\otimes W)\otimes W^{\prime}=\mathcal{D}\otimes(W\otimes_{\infty
}W^{\prime})
\]
for any Weil algebra $W^{\prime}$, so that we have the desired result.
\end{proof}

\section{Microlinearity\label{3}}

A Fr\"{o}licher space $X$ is called \textit{microlinear} providing that any
finite limit diagram $\mathcal{D}$ in $\mathbf{W}$ yields a limit diagram
$X\otimes\mathcal{D}$ in $\mathbf{FS}$, where $X\otimes\mathcal{D}$ is
obtained from $\mathcal{D}$ by putting $X\otimes$ to the left of every object
and every morphism in $\mathcal{D}$.

The following result should be obvious.

\begin{proposition}
\label{t3.1}Convenient vector spaces are microlinear.
\end{proposition}

\begin{corollary}
$C^{\infty}$-manifolds in the sense of \cite{kri} (cf. Section 27) are microlinear.
\end{corollary}

\begin{proposition}
\label{t3.2}If $X$ is a Weil exponentiable and microlinear Fr\"{o}licher
space, then so is $X\otimes W$ for any Weil algebra $W$.
\end{proposition}

\begin{proof}
This follows simply from Proposition \ref{t2.2.1} and Lemma \ref{t2.3.2}.
\end{proof}

\begin{proposition}
\label{t3.3}If $X$ and $Y$ are microlinear Fr\"{o}licher spaces, then so is
$X\times Y$.
\end{proposition}

\begin{proof}
This follows simply from Theorem \ref{t2.1} and the familiar fact that the
functor $\cdot\times\cdot:$ $\mathbf{FS}\times\mathbf{FS}\rightarrow
\mathbf{FS}$ preserves limits.
\end{proof}

\begin{proposition}
\label{t3.4}If $X$ is a Weil exponentiable and microlinear Fr\"{o}licher
space, then so is $X^{Y}$ for any Fr\"{o}licher space $Y$.
\end{proposition}

\begin{proof}
This follows simply from (\ref{2.3}), Proposition \ref{t2.2.3} and Lemma
\ref{t2.3.1}.
\end{proof}

We recapitulate:

\begin{theorem}
\label{t3.5}Weil exponentiable and microlinear Fr\"{o}licher spaces, together
with smooth mappings among them, form a cartesian closed subcategory
$\mathbf{FS}_{\mathbf{WE,ML}}$\ of the category $\mathbf{FS}$.
\end{theorem}

We note in passing that microlinearity is closed under transversal limits.

\begin{theorem}
\label{t3.6}If the limit of a diagram $\mathcal{F}$ of microlinear
Fr\"{o}licher spaces is transversal, then it is microlinear.
\end{theorem}

\begin{proof}
Let $X$ be the limit of the diagram $\mathcal{F}$, i.e.,
\[
X=\mathrm{Lim\ }\mathcal{F}%
\]
Let $W$ be the limit of an arbitrarily given finite diagram $\mathcal{D}$ of
Weil algebras, i.e.,
\[
W=\mathrm{Lim\ }\mathcal{D}%
\]
We denote by $\mathcal{F}\otimes\mathcal{D}$ the diagram obtained from the
diagrams $\mathcal{F}$ and $\mathcal{D}$ by the application of the bifunctor
$\otimes:\mathbf{FS}\times\mathbf{W\rightarrow FS}$. By recalling that double
limits in a complete category commute (cf. Section 2 of Chapter IX of
\cite{mac}), we have
\begin{align*}
&  \mathrm{Lim\ }(X\otimes\mathcal{D})\\
&  =\mathrm{Lim\ }((\mathrm{Lim\ }\mathcal{F})\otimes\mathcal{D)}\\
&  =\mathrm{Lim}_{\mathcal{D}}\mathrm{\ Lim}_{\mathcal{F}}\mathrm{\ }%
(\mathcal{F}\otimes\mathcal{D}\mathbb{)}\\
&  \text{[since }\mathrm{Lim\ }\mathcal{F}\text{ is the transversal limit]}\\
&  =\mathrm{Lim}_{\mathcal{F}}\mathrm{\ Lim}_{\mathcal{D}}\mathrm{\ }%
(\mathcal{F}\otimes\mathcal{D}\mathbb{)}\\
&  \text{[since double limits commute]}\\
&  =\mathrm{Lim\ }(\mathcal{F}\otimes(\mathrm{Lim\ }\mathcal{D}))\\
&  \text{[since every object in }\mathcal{F}\text{ is microlinear]}\\
&  =\mathrm{Lim\ }(\mathcal{F}\otimes W)\\
&  =(\mathrm{Lim\ }\mathcal{F)}\otimes W\\
&  \text{[since }\mathrm{Lim\ }\mathcal{F}\text{ is the transversal limit]}\\
&  =X\otimes W
\end{align*}
Therefore we have the desired result.
\end{proof}

\begin{proposition}
\label{t3.7}If a Weil exponentiable Fr\"{o}licher space $X$ is microlinear,
then any finite limit diagram $\mathcal{D}$ in $\mathbf{W}$ yields a
transversal limit diagram $X\otimes\mathcal{D}$ in $\mathbf{FS}$.
\end{proposition}

\begin{proof}
By the same token as in the proof of Lemma \ref{t2.3.2}.
\end{proof}

\section{The Embedding into the Cahiers Topos\label{4}}

Let $\mathbf{D}$ be the full subcategory of the category of $\mathcal{C}%
^{\infty}$-algebras in form $\mathcal{C}^{\infty}(\mathbb{R}^{n}%
)\otimes_{\infty}W$ with a natural number $n$ and a Weil algebra $W$. Now we
would like to extend the Weil prolongation $\mathbf{FS}_{\mathbf{WE}}%
\times\mathbf{W}\overset{\otimes}{\mathbf{\rightarrow}}\mathbf{FS}%
_{\mathbf{WE}}$ to a bifunctor $\mathbf{FS}_{\mathbf{WE}}\times\mathbf{D}%
\overset{\otimes}{\mathbf{\rightarrow}}\mathbf{FS}_{\mathbf{WE}}$. On objects
we define
\begin{equation}
X\otimes C=X^{\mathbb{R}^{n}}\otimes W\label{4.1}%
\end{equation}
for any Weil exponentiable Fr\"{o}licher space $X$ and any $C=\mathcal{C}%
^{\infty}(\mathbb{R}^{n})\otimes_{\infty}W$. By Proposition \ref{t2.2.3}
$X^{\mathbb{R}^{n}}$ is Weil exponentiable, so that $X\otimes C$ is Weil
exponentiable by Proposition \ref{t2.2.1}. It is easy to see that the right
hand of (\ref{4.1}) is functorial in $X$, but we have not so far succeeded in
establishing its functoriality in $C$. Therefore we pose it as a conjecture.

\begin{conjecture}
\label{t4.1}The right hand of (\ref{4.1}) is functorial in $C$, so that we
have a bifunctor $\mathbf{FS}_{\mathbf{WE}}\times\mathbf{D}\overset{\otimes
}{\mathbf{\rightarrow}}\mathbf{FS}_{\mathbf{WE}}$.
\end{conjecture}

In the following we will assume that the conjecture is really true. We define
the functor $\mathbf{J:FS}_{\mathbf{WE}}\rightarrow\mathbf{Sets}^{\mathbf{D}}$
to be the exponential adjoint to the composite
\[
\mathbf{FS}_{\mathbf{WE}}\times\mathbf{D}\overset{\otimes}{\mathbf{\rightarrow
}}\mathbf{FS}_{\mathbf{WE}}\rightarrow\mathbf{Sets}%
\]
where $\mathbf{FS}_{\mathbf{WE}}\rightarrow\mathbf{Sets}$ is the
underlying-set functor. Now we have

\begin{proposition}
\label{t4.2}For any Weil-exponentiable Fr\"{o}licher space $X$ and any object
$\mathcal{C}^{\infty}(\mathbb{R}^{n})\otimes_{\infty}W$ in $\mathbf{D} $, we
have
\[
\mathbf{J}(X)^{\hom_{\mathbf{D}}(\mathcal{C}^{\infty}(\mathbb{R}^{n}%
)\otimes_{\infty}W,\cdot)}=\mathbf{J}(X\otimes(\mathcal{C}^{\infty}%
(\mathbb{R}^{n})\otimes_{\infty}W))
\]

\end{proposition}

\begin{proof}
For any object $\mathcal{C}^{\infty}(\mathbb{R}^{m})\otimes_{\infty}W^{\prime
}$ in $\mathbf{D}$, we have
\begin{align*}
&  \mathbf{J}(X)^{\hom_{\mathbf{D}}(\mathcal{C}^{\infty}(\mathbb{R}%
^{n})\otimes_{\infty}W,\cdot)}(\mathcal{C}^{\infty}(\mathbb{R}^{m}%
)\otimes_{\infty}W^{\prime})\\
&  =\hom_{\mathbf{sets}^{\mathbf{D}}}(\hom_{\mathbf{D}}(\mathcal{C}^{\infty
}(\mathbb{R}^{m})\otimes_{\infty}W^{\prime},\cdot),\mathbf{J}(X)^{\hom
_{\mathbf{D}}(\mathcal{C}^{\infty}(\mathbb{R}^{n})\otimes_{\infty}W,\cdot)})\\
&  \text{[By Yoneda Lemma]}\\
&  =\hom_{\mathbf{sets}^{\mathbf{D}}}(\hom_{\mathbf{D}}(\mathcal{C}^{\infty
}(\mathbb{R}^{m})\otimes_{\infty}W^{\prime},\cdot)\times\hom_{\mathbf{D}%
}(\mathcal{C}^{\infty}(\mathbb{R}^{n})\otimes_{\infty}W,\cdot),\mathbf{J}%
(X))\\
&  =\hom_{\mathbf{sets}^{\mathbf{D}}}(\hom_{\mathbf{D}}(\mathcal{C}^{\infty
}(\mathbb{R}^{m})\otimes_{\infty}\mathcal{C}^{\infty}(\mathbb{R}^{n}%
)\otimes_{\infty}W\otimes_{\infty}W^{\prime},\cdot),\mathbf{J}(X))\\
&  =\hom_{\mathbf{sets}^{\mathbf{D}}}(\hom_{\mathbf{D}}(\mathcal{C}^{\infty
}(\mathbb{R}^{m}\times\mathbb{R}^{n})\otimes_{\infty}W\otimes_{\infty
}W^{\prime},\cdot),\mathbf{J}(X))\\
&  =\mathbf{J}(X)(\mathcal{C}^{\infty}(\mathbb{R}^{m}\times\mathbb{R}%
^{n})\otimes_{\infty}W\otimes_{\infty}W^{\prime})\\
&  \text{[By Yoneda Lemma]}\\
&  =X^{\mathbb{R}^{m}\times\mathbb{R}^{n}}\otimes(W\otimes_{\infty}W^{\prime
})\\
&  =(X^{\mathbb{R}^{n}}\otimes W)^{\mathbb{R}^{m}}\otimes W^{\prime}\\
&  =(X\otimes(\mathcal{C}^{\infty}(\mathbb{R}^{n})\otimes_{\infty
}W))^{\mathbb{R}^{m}}\otimes W^{\prime}\\
&  =\mathbf{J}(X\otimes(\mathcal{C}^{\infty}(\mathbb{R}^{n})\otimes_{\infty
}W))(\mathcal{C}^{\infty}(\mathbb{R}^{m})\otimes_{\infty}W^{\prime})
\end{align*}
Therefore we have the desired result.
\end{proof}

\begin{proposition}
\label{t4.3}For any Weil-exponentiable Fr\"{o}licher spaces $X$ and $Y$, we
have
\[
\mathbf{J}(X\times Y)=\mathbf{J}(X)\times\mathbf{J}(Y)
\]

\end{proposition}

\begin{proof}
For any object $\mathcal{C}^{\infty}(\mathbb{R}^{n})\otimes_{\infty}W$ in
$\mathbf{D}$, we have
\begin{align*}
&  \mathbf{J}(X\times Y)(\mathcal{C}^{\infty}(\mathbb{R}^{n})\otimes_{\infty
}W)\\
&  =(X\times Y)^{\mathbb{R}^{n}}\otimes W\\
&  =(X^{\mathbb{R}^{n}}\times Y^{\mathbb{R}^{n}})\otimes W\\
&  =(X^{\mathbb{R}^{n}}\otimes W)\times(Y^{\mathbb{R}^{n}}\otimes W)\\
&  \text{[By Theorem \ref{t2.1}]}\\
&  =\mathbf{J}(X)(\mathcal{C}^{\infty}(\mathbb{R}^{n})\otimes W)\times
\mathbf{J}(Y)(\mathcal{C}^{\infty}(\mathbb{R}^{n})\otimes_{\infty}W)\\
&  =(\mathbf{J}(X)\times\mathbf{J}(Y))(\mathcal{C}^{\infty}(\mathbb{R}%
^{n})\otimes_{\infty}W)
\end{align*}
Therefore we have the desired result.
\end{proof}

\begin{proposition}
\label{t4.4}For any Weil exponentiable Fr\"{o}licher spaces $X$ and $Y$, we
have the following isomorphism in $\mathbf{Sets}^{\mathbf{D}}$:
\[
\mathbf{J}(X^{Y})=\mathbf{J}(X)^{\mathbf{J}(Y)}%
\]

\end{proposition}

\begin{proof}
Let $\mathcal{C}^{\infty}(\mathbb{R}^{n})\otimes_{\infty}W$ be an object in
$\mathbf{D}$. On the one hand, we have
\begin{align*}
&  \mathbf{J}(X^{Y})(\mathcal{C}^{\infty}(\mathbb{R}^{n})\otimes_{\infty}W)\\
&  =(X^{Y})^{\mathbb{R}^{n}}\otimes W\\
&  =X^{Y\times\mathbb{R}^{n}}\otimes W\\
&  =(X^{\mathbb{R}^{n}}\otimes W)^{Y}%
\end{align*}
On the other hand, we have
\begin{align*}
&  \mathbf{J}(X)^{\mathbf{J}(Y)}(\mathcal{C}^{\infty}(\mathbb{R}^{n}%
)\otimes_{\infty}W)\\
&  =\hom_{\mathbf{Sets}^{\mathbf{D}}}(\hom_{\mathbf{D}}(\mathcal{C}^{\infty
}(\mathbb{R}^{n})\otimes_{\infty}W,\cdot),\mathbf{J}(X)^{\mathbf{J}(Y)})\\
&  \text{[By Yoneda Lemma]}\\
&  =\hom_{\mathbf{Sets}^{\mathbf{D}}}(\hom_{\mathbf{D}}(\mathcal{C}^{\infty
}(\mathbb{R}^{n})\otimes_{\infty}W,\cdot)\times\mathbf{J}(Y),\mathbf{J}(X))\\
&  =\hom_{\mathbf{Sets}^{\mathbf{D}}}(\mathbf{J}(Y),\mathbf{J}(X)^{\hom
_{\mathbf{D}}(\mathcal{C}^{\infty}(\mathbb{R}^{n})\otimes_{\infty}W,\cdot)})\\
&  =\hom_{\mathbf{Sets}^{\mathbf{D}}}(\mathbf{J}(Y),\mathbf{J}(X\otimes
(\mathcal{C}^{\infty}(\mathbb{R}^{n})\otimes_{\infty}W)))\\
&  \text{[By Proposition \ref{t4.2}]}\\
&  =\hom_{\mathbf{FS}_{\mathbf{WE}}}(Y,X\otimes(\mathcal{C}^{\infty
}(\mathbb{R}^{n})\otimes_{\infty}W))\\
&  \text{[since }\mathbf{J}\text{\ is full and faithful]}%
\end{align*}
Therefore we have the desired result.
\end{proof}

\begin{theorem}
\label{t4.5}The functor $\mathbf{J}:\mathbf{FS}_{\mathbf{WE}}\rightarrow
\mathbf{Sets}^{\mathbf{D}}$ preserves the cartesian closed structure. In other
words, it preserves finite products and exponentials. It is full and faithful.
It sends the Weil prolongation to the exponentiation by the corresponding
infinitesimal object.
\end{theorem}

\begin{proof}
The first statement follows from Propositions \ref{t4.3} and \ref{t4.4}. The
second statement that it is full and faithful follows by the same token as in
\cite{kock1} and \cite{kock2}. The final statement follows from Proposition
\ref{t4.2}.
\end{proof}

Now we are concerned with microlinearity. First we will establish

\begin{proposition}
\label{t4.6}The functor $\mathbf{J}:\mathbf{FS}_{\mathbf{WE}}\rightarrow
\mathbf{Sets}^{\mathbf{D}}$ preserves transversal limit diagrams. In other
words, the functor $\mathbf{J}$ always sends a transversal limit diagram of
Fr\"{o}licher spaces lying in $\mathbf{FS}_{\mathbf{WE}}$ to the limit diagram
in $\mathbf{Sets}^{\mathbf{D}}$.
\end{proposition}

\begin{proof}
We should show that a transversal limit diagram $\mathcal{D}$ in
$\mathbf{FS}_{\mathbf{WE}}$\ always yields a limit diagram $\mathbf{J}%
(\mathcal{D})$ in $\mathbf{Sets}^{\mathbf{D}}$. To this end, it suffices to
show (cf. \cite{macmoe}, pp. 22-23) that $\mathbf{J}(\mathcal{D})(C)$ is a
limit diagram in $\mathbf{Sets}$\ for any object $C=\mathcal{C}^{\infty
}(\mathbb{R}^{n})\otimes_{\infty}W$ in $\mathbf{D}$. Since the forgetful
functor $\mathbf{FS\rightarrow Sets}$ preserves limits, we have only to note
that $\mathcal{D}^{\mathbb{R}^{n}}\otimes W$ is a limit diagram in
$\mathbf{FS}$, which follows readily from Lemma \ref{t2.3.1}.
\end{proof}

\begin{theorem}
\label{t4.7}The functor $\mathbf{J}:\mathbf{FS}_{\mathbf{WE}}\rightarrow
\mathbf{Sets}^{\mathbf{D}}$ preserves microlinearity.
\end{theorem}

\begin{proof}
This follows simply Propositions \ref{t3.7}, \ref{t4.2} and \ref{t4.6}.
\end{proof}

The site of definition for the Cahiers topos $\mathcal{C}$ is the dual
category $\mathbf{D}^{\mathrm{op}}$\textbf{\ }of the category $\mathbf{D}$
together with the open-cover topology, so that we have the canonical
embedding
\[
\mathcal{C\hookrightarrow}\mathbf{Sets}^{\mathbf{D}}%
\]
By the same token as in \cite{kock1} and \cite{kock2} we can see that the
functor $\mathbf{J}:\mathbf{FS}_{\mathbf{WE}}\rightarrow\mathbf{Sets}%
^{\mathbf{D}}$ factors in the above embedding. The resulting functor is
denoted by $\mathbf{J}_{\mathcal{C}}$. Since the above embedding creates
limits and exponentials, Theorems \ref{t4.5} and \ref{t4.7} yields directly

\begin{theorem}
\label{t4.8}The functor $\mathbf{J}_{\mathcal{C}}:\mathbf{FS}_{\mathbf{WE}%
}\rightarrow\mathcal{C}$ preserves the cartesian closed structure. In other
words, it preserves finite products and exponentials. It is full and faithful.
It sends the Weil prolongation to the exponentiation by the corresponding
infinitesimal object. The functor $\mathbf{J}_{\mathcal{C}}$ preserves
transversal limit diagrams and microlinearity.
\end{theorem}

\end{document}